\documentclass{article}
\usepackage{latexsym}
\usepackage{amsmath}
\usepackage{amssymb}
\usepackage{theorem}
\makeatletter
\@ifpackageloaded{amsgen}{}{\input{amsgen.sty}}
\renewcommand{\theoremstyle}[1]{%
  \@ifundefined{th@#1}{%
    \PackageWarning{amsthm}{Unknown theoremstyle `#1'}%
    \theorem@style{plain}%
  }{%
    \theorem@style{#1}%
  }%
}
\newtoks\theorem@style
\theorem@style{plain}
\newtoks\theorembodyfont
\theorembodyfont{\itshape}
\newtoks\theoremheadfont
\theoremheadfont{\bfseries}
\newtoks\theoremnotefont
\theoremnotefont{\bfseries}
\newtoks\theoremheadpunct
\theoremheadpunct{.}
\newskip\theorempreskipamount \theorempreskipamount\topsep
\newskip\theorempostskipamount \theorempostskipamount\topsep
\renewcommand{\newtheorem}{\@ifstar{\@xnthm *}{\@xnthm \relax}}
\def\@xnthm#1#2{%
  \let\@tempa\relax
  \@xp\@ifdefinable\csname #2\endcsname{%
    \global\@xp\let\csname end#2\endcsname\@endtheorem
    \ifx *#1
      \edef\@tempa##1{%
        \gdef\@xp\@nx\csname#2\endcsname{%
          \@nx\@thm{\@xp\@nx\csname th@\the\theorem@style\endcsname}%
            {}{##1}}}%
    \else 
      \def\@tempa{\@oparg{\@ynthm{#2}}[]}%
    \fi
  }%
  \@tempa
}
\def\@ynthm#1[#2]#3{%
  \ifx\relax#2\relax
    \def\@tempa{\@oparg{\@xthm{#1}{#3}}[]}%
  \else
    \@ifundefined{c@#2}{%
      \def\@tempa{\@nocounterr{#2}}%
    }{%
      \@xp\xdef\csname the#1\endcsname{\@xp\@nx\csname the#2\endcsname}%
      \toks@{#3}%
      \@xp\xdef\csname#1\endcsname{%
        \@nx\@thm{%
          \let\@nx\thm@swap
            \if S\thm@swap\@nx\@firstoftwo\else\@nx\@gobble\fi
          \@xp\@nx\csname th@\the\theorem@style\endcsname}%
            {#2}{\the\toks@}}%
      \let\@tempa\relax
    }%
  \fi
  \@tempa
}
\def\@xthm#1#2[#3]{%
  \ifx\relax#3\relax
    \newcounter{#1}%
  \else
    \newcounter{#1}[#3]%
    \@xp\xdef\csname the#1\endcsname{\@xp\@nx\csname the#3\endcsname
      \@thmcountersep\@thmcounter{#1}}%
  \fi
  \toks@{#2}%
  \@xp\xdef\csname#1\endcsname{%
    \@nx\@thm{%
      \let\@nx\thm@swap
        \if S\thm@swap\@nx\@firstoftwo\else\@nx\@gobble\fi
      \@xp\@nx\csname th@\the\theorem@style\endcsname}%
      {#1}{\the\toks@}}%
}
\def\@thm#1#2#3{\normalfont
  \trivlist
  \labelsep.5em\relax \let\thmheadnl\relax
  \let\theoremindent\noindent 
  \let\thm@swap\@gobble
  \theoremheadfont{\bfseries}
  \theoremheadpunct{.}
  \theorempreskipamount\topsep
  \theorempostskipamount\theorempreskipamount
  #1
  \@topsep \theorempreskipamount               
  \@topsepadd \theorempostskipamount           
  \def\@tempa{#2}\ifx\@empty\@tempa
    \def\@tempa{\@oparg{\@begintheorem{#3}{}}[]}%
  \else
    \refstepcounter{#2}%
    \def\@tempa{\@oparg{\@begintheorem{#3}{\csname the#2\endcsname}}[]}%
  \fi
  \@tempa
}
\let\@ythm\relax
\let\thmname\@iden \let\thmnumber\@iden \let\thmnote\@iden
\def\thmhead@plain#1#2#3{%
  \thmname{#1}\thmnumber{ #2}\thmnote{ {\the\theoremnotefont(#3)}}}
\let\thmhead\thmhead@plain
\def\swappedhead#1#2#3{%
  \thmnumber{#2}\thmname{. #1}\thmnote{ {\the\theoremnotefont(#3)}}}
\let\thmheadnl\relax
\def\@begintheorem#1#2[#3]{%
  \item[\normalfont 
  \hskip\labelsep
  \the\theoremheadfont
  \theoremindent
  \@ifempty{#1}{\let\thmname\@gobble}{\let\thmname\@iden}%
  \@ifempty{#2}{\let\thmnumber\@gobble}{\let\thmnumber\@iden}%
  \@ifempty{#3}{\let\thmnote\@gobble}{\let\thmnote\@iden}%
  \thm@swap\swappedhead\thmhead{#1}{#2}{#3}%
  \the\theoremheadpunct]%
  \thmheadnl 
  \ignorespaces}
\def\nonslanted{\relax
  \@xp\let\@xp\@tempa\csname\f@shape shape\endcsname
  \ifx\@tempa\itshape\upshape
  \else\ifx\@tempa\slshape\upshape\fi\fi}
\def\swapnumbers{\edef\thm@swap{\if S\thm@swap N\else S\fi}}
\def\thm@swap{N}%
\let\@opargbegintheorem\relax
\def\th@plain{%
  \itshape 
}
\def\th@definition{%
  \normalfont 
}
\def\th@remark{%
  \theoremheadfont{\itshape}%
  \normalfont 
  \theorempreskipamount\topsep
  \divide\theorempreskipamount\tw@
  \theorempostskipamount\theorempreskipamount
}
\def\@endtheorem{\endtrivlist\@endpefalse }
\newcommand{\newtheoremstyle}[9]{%
  \@ifempty{#5}{\dimen@\z@skip}{\dimen@#5\relax}%
  \ifdim\dimen@=\z@
    \toks@{#4\let\theoremindent\noindent}%
  \else
    \toks@{#4\def\theoremindent{\noindent\hbox to#5{}}}%
  \fi
  \def\@tempa{#8}\ifx\space\@tempa
    \toks@\@xp{\the\toks@ \labelsep\fontdimen\tw@\font\relax}%
  \else
    \def\@tempb{\newline}%
    \ifx\@tempb\@tempa
      \toks@\@xp{\the\toks@ \labelsep\z@skip
        \def\thmheadnl{\leavevmode\newline}}%
    \else
      \toks@\@xp{\the\toks@ \labelsep#8\relax}%
    \fi
  \fi
  \@temptokena{%
    \theorempreskipamount#2\relax
    \theorempostskipamount#3\relax
    \theoremheadfont{#6}\theoremheadpunct{#7}%
  }%
  \@ifempty{#9}{%
    \let\thmhead\thmhead@plain
  }{%
    \@namedef{thmhead@#1}##1##2##3{#9}%
    \@temptokena\@xp{\the\@temptokena
      \@xp\let\@xp\thmhead\csname thmhead@#1\endcsname}%
  }%
  \@xp\xdef\csname th@#1\endcsname{\the\toks@ \the\@temptokena}%
}
\theoremstyle{remark}
\newtheorem{remark}{Remark}
\newtheorem{proof}[remark]{Proof}

\theoremstyle{plain}
\newtheorem{lemma}{Lemma}
\newtheorem{thm}[lemma]{Theorem}

\newtheorem{proposition}[lemma]{Proposition}
\newtheorem{corollary}[lemma]{Corollary}
\newtheorem{claim}[lemma]{Claim}

\theoremstyle{definition}

\makeatletter
\renewcommand{\subsection}{\@startsection{subsection}{2}{\z@}{-3.5ex}{-0ex}{
\bf}}
\makeatother
\def\sss{\subsection{}}

\def\wt{\mathrm{wt}}

\def\bfm{{\bf m}}\def\bfn{{\bf n}}\def\bfmp{{\bf m'}}
\def\bfnw{{\bf n}_{\tilde w}}\def\bfp{{\bf p}}\def\bfr{{\bf r}}
\def\zpn{\Z_{\geq 0}^N}\def\zprp{\Z_{\geq 0}^{\rp}}

\def\Hom{\mathrm{Hom}}

\def\B{{\cal B}}

\def\ga{\Gamma}\def\gaf{\overrightarrow{\ga}}

\def\uqnm{U(\nm)}

\def\uq{U_q(\g)}
\def\uqn{U_q(\n)}
\def\uqnm{U_q(\nm)}
\def\uo{U_q^0}

\def\uqb{U_q(\b)}
\def\uqbm{U_q(\b^-)}

\def\lq{V_q}
\def\lql{\lq(\lambda)}

\def\dim{\hbox{dim}}\def\Hom{\hbox{Hom}}\def\Ext{\hbox{Ext}}
\def\tr{\mathrm{tr}}

\def\N{{\mathbb N}}\def\Z{{\mathbb Z}}%
\def\Q{{\mathbb Q}}%
\def\n{{\mathfrak n}}%
\def\C{{\mathbb C}}%
\def\g{{\mathfrak g}}%
\def\h{{\mathfrak h}}%
\def\b{{\mathfrak b}}\def\nm{{\mathfrak n}^-}
\def\lat{{\cal L}}

\def\rp{R^+}

\begin{document}
\title{A multiplicative property of quantum flag minors}
\author{Ph. Caldero\thanks{Supported in part by the EC TMR network
 ``Algebraic Lie Representations", contract no. ERB FMTX-CT97-0100}
}
\maketitle
\begin{abstract} We study the multiplicative properties of the quantum dual canonical basis ${\cal B}^*$ associated to a semi-simple complex Lie group $G$. We provide a subset $D$ of ${\cal B}^*$ such that the following property holds : if two elements $b$, $b'$ in ${\cal B}^*$ $q$-commute and if one of these elements is in $D$, then the product $bb'$ is in ${\cal B}^*$ up to a power of $q$, where $q$ the quantum parameter. If $G$ is SL$_n$, then $D$ is the set of so-called quantum flag minors and we obtain a generalization of a result of Leclerc-Nazarov-Thibon, [11].
\end{abstract}
\setcounter{section}{-1}
\section{Introduction}
\sss Let $G$ be a semisimple complex Lie group and fix a maximal unipotent
subgroup $U^-$ of $G$. Let $\g$ and $\nm$ be respectively the Lie algebras of
$G$ and $U^-$. G. Lusztig and M. Kashiwara have constructed the so-called
canonical basis $\B$ of the enveloping algebra $U(\nm)$ of $\nm$, which has
 properties of compatibility with standard filtrations.\par
Let $\C[U^-]$ be the $\C$-algebra of regular functions on $U^-$. Then, the
action of $U^-$ on itself by left multiplication provides an action of
$U(\nm)$ on $\C[U^-]$ by differential operators. Now, consider the pairing
$U(\nm)\times \C[U^-]\rightarrow\C$, $\delta\times f\mapsto
\delta(f)(\hbox{e})$, where $e$ is the identity of $U^-$. Then, this pairing
provides the so-called dual canonical basis $\B^*$ of $\C[U^-]$. This
article is concerned with some multiplicative properties of this basis.
\sss Let $q$ be an indeterminate and let $U_q(\nm)$, $\C_q[U^-]$ be
respectively the quantum analogue of the classical objects $U(\nm)$ and
$\C[U^-]$. We still note $\B$, resp. $\B^*$, the canonical basis, resp. the
dual canonical basis, of $U_q(\nm)$, resp. $\C_q[U^-]$. We say that two
elements $b$ and $b'$ of $\B^*$ $q$-commute if $bb'=q^mb'b$, for an
integer $m$. We say that they are multiplicative if $q^nbb'$ belongs to
$\B^*$ for an integer $n$. It is known, see [19], that if two elements are multiplicative, then they $q$-commute. The converse was believed to be true until Bernard Leclerc found counter examples. \par
Suppose that two elements $b$ and $b'$ of the dual canonical basis $q$-commute. Let's discuss now in which cases they are known to be multiplicative.\par\noindent
1) for all $b$, $b'$, if $\g$ is of type A$_n$, $n\leq 3$, B$_2$, [1], see also
[4].\par\noindent
2) if $\g$ is of type A$_n$ and $b$ is a small quantum minor,
[19].\par\noindent
3) if $\g$ is of type A$_n$ and $b$, $b'$ are quantum flag
minors, [11].\par
\sss Let's present the results of this article. Let $W$ be the Weyl group of
$\g$ and let $w_0$ be its longest element. For each reduced decomposition
$\tilde w_0$ of $w_0$, we have constructed in [4] a subalgebra
$A_{\tilde w_0}$ of $\C_q[U^-]$ such that:\par\noindent
1) The space $A_{\tilde w_0}$ is generated by a part of $\B^*$,\par\noindent
2) pair of elements in $\B^*\cap A_{\tilde w_0}$ are multiplicative,\par\noindent
3) the algebras $A_{\tilde w_0}$ and $\C_q[U^-]$ are equal up to
localization.\par\noindent
They are called adapted algebras associated to a reduced
decomposition of $w_0$. They are connected to the more general theory of
cluster algebras, [8], see 5.2.\par
The main result of the article is the following:
\begin{thm} Let $\tilde w_0$ be a reduced decomposition corresponding to an
orientation of the Coxeter graph of $\g$. Let $b$ and $b'$ be two
$q$-commuting elements of the dual canonical basis. If $b$ is in $A_{\tilde
w_0}$ then $b$ and $b'$ are multiplicative.
\end{thm}
As a particular case, we obtain:
\begin{corollary} Let $\g$ be of type $A_n$. Let $b$ be a quantum flag minor
and let $b'$ be any element of $\B^*$ which $q$-commutes with $b$, then $b$
and $b'$ are multiplicative.
\end{corollary}
We refer to [11] for motivations of this result. Actually, it provides a criterion of irreducibility for module on the affine Hecke algebra of type A which are induced by the so-called evaluation modules.
\par\noindent
Let's sketch the proof of the theorem. If two elements $b$ and $b'$ of the
dual canonical basis $q$-commute, then, the only property we have to obtain
in order to prove that $b$ and $b'$ are multiplicative is
$$q^nbb'\in b''+q{\cal L}^*,\leqno (0.3.1)$$
where $n$ is an integer, $b''$ is in $\B^*$ and where ${\cal L}^*$ is the
$\Z[q]$-lattice generated by $\B^*$.
The natural question is: how to control the powers of $q$ in the
multiplications of elements of the dual canonical basis ? The control of
these powers are based on two main ideas:\par\noindent
1) Kashiwara proved that bases of integrable modules of the quantized
enveloping algebra $\uq$ cristallizes at $q=0$, with compatibility with the
tensor product. To be more precise, let $P^+$ be the semigroup of integral
dominant weights and let $\overline b$, $\overline b'$ be the corresponding
elements in the crystal bases $\B(\lambda)$, $\B(\lambda')$,
of the integrable modules of highest weight, respectively, $\lambda$ and
$\lambda'$ in $P^+$. We can  assert that if  $\overline
b\otimes\overline b'$ belongs to the connected component of the crystal
$\B(\lambda)\otimes\B(\lambda')$ corresponding to $\B(\lambda+\lambda')$,
then (0.3.1) holds. The property $\overline
b\otimes\overline b'\in\B(\lambda+\lambda')$ can be checked easily via the Littelmann's
path model of the crystal basis, [13], [14], by comparing chains of elements of the Weyl
group for the Bruhat ordering.\par\noindent
2) The quiver approach of the algebra $\C_q[U^-]$ enables to interpret powers
of $q$ which appear in multiplications in terms of dim Hom$(M,N)$
and dim Ext$^1(M,N)$, where $M$ and $N$ are representations of a quiver. An
important tool is that the map  dim Hom$(M,?)$ is increasing for the
so-called degeneration ordering, see [3].
 \section{Notations.}
\sss Let $\g$ be a semi-simple Lie $\C$-algebra of rank $n$ with Cartan
matrix $A=(a_{ij})$.
We fix a Cartan subalgebra $\h$ of $\g$. Let $\g=\nm \oplus\h \oplus\n$
be the triangular decomposition and set $N:=\dim\n$. Let $\{\alpha_i\}_i$ be
a basis of the
root system
$R$ resulting from this decomposition and let $\rp$ be the set of
positive
roots.
Let $P$ be the weight $\Z$-lattice generated by the fundamental weights
$\varpi_i$, $i\in I:=\Z\cap[1,n]$ and set $P^+:=\sum_i\,\N\varpi_i$.
Let $W$ be the Weyl group, generated by the
reflections corresponding to the simple roots $s_i:=s_{\alpha_i}$, with
longest element
$w_0$.
We note $<\,,\,>$ the $W$-invariant form on $P$.
\sss In this section we define the quantized enveloping algebra of $\g$ and
the properties of its Poincar\'e-Birkhoff-Witt basis. We refer to [6] for
precisions and proofs.\par\noindent
Let $q$ be an indeterminate. Let $U_q(\g)$ be the quantized enveloping Hopf
$\Q(q)$-algebra as defined in [6]. Let $\uqn$, resp $\uqnm$, be the
upper, resp
lower, ``nilpotent" subalgebra of $\uq$. The algebra $\uqn$, resp $\uqnm$,
is generated by $E_i$, resp $F_i$, $1\leq i\leq n$ with quantum Serre
relations. For all $\lambda$
in $Q:=\oplus_i\Z\alpha_i$, let $K_\lambda$ be the corresponding element in
the algebra $\uo=\Q(q)[Q]$ of the torus of $\uq$.
\par
For all $\mu\in Q$, let $\uqn_\mu$ be the subspace of $\uqn$ generated by
the products $E_{i_1}^{n_1}\ldots E_{i_k}^{n_k}$ such that
$\sum_ln_l\alpha_l=\mu$. An element $X$ of $\uqn_\mu$ will be called
(homogeneous) element of weight $\mu$. We set $\wt(X):=\mu$.\par
Recall the triangular decomposition
$\uq=\uqnm\otimes\uo\otimes\uqn$. We define the following subalgebras of
$\uq$:
$$\uqb=\uqn\otimes\uo, \hskip 15mm \uqbm=\uqnm\otimes\uo.$$
As in [22], [16], we introduce Lusztig's automorphisms $T_i$,
$1\leq
i\leq
n$, which define a braid action on $\uq$ by
$$T_i(E_i)=-F_iK_{\alpha_i},\hskip 5mm T_i(E_j)=\sum_{s=0}^{-a_{ij}}
(-1)^{-a_{ij}-s}q_{\alpha_i}^{a_{ij}+s}E_i^{(s)}E_jE_i^{(-a_{ij}-s)},\;
1\leq i,j\leq n,\, i\not=j\leqno (1.2.1)$$
$$T_i(F_i)=-K_{-\alpha_i}E_i,\hskip 3mm T_i(F_j)=(-1)^{-a_{ij}-s}
\sum_{s=0}^{-a_{ij}}
q_{\alpha_i}^{-a_{ij}-s}F_i^{(-a_{ij}-s)}F_jF_i^{(s)},\; 1\leq i,j\leq n,\,
i\not=j\leqno (1.2.2)$$
$$T_i(K_{\alpha_j})=K_{s_i(\alpha_j)},\,1\leq i,j\leq n,\leqno (1.2.3)$$
where $E_{\alpha}^{(k)}={1\over[k]_{q_{\alpha}}!}E_{\alpha}^{k}$,
$[k]_{q_{\alpha}}!=[k]_{q_{\alpha}}
[k-1]_{q_{\alpha}}\ldots[1]_{q_{\alpha}}$,\par\noindent
 $[k]_{q_{\alpha}}={q_{\alpha}^n-q_{\alpha}^{-n}\over
q_{\alpha}-q_{\alpha}^{-1}}$,
 $q_{\alpha}=q^{{(\alpha,\alpha)\over 2}}$.
Fix a reduced decomposition $\tilde w_0=s_{i_1}s_{i_2}\ldots s_{i_N}$ of
$w_0$. For each $k$, $1\leq k\leq N$, set $\beta_k:=s_{i_1}\ldots
s_{i_{k-1}}(\alpha_k)$. It is well known, [17], that the $\{\beta_k,\,
1\leq k\leq N\}$ is the set of positive roots and that
$$\beta_1<\beta_2<\ldots<\beta_N$$
defines a so-called convex ordering on $\rp$. This ordering will identify
the semigroup $\zprp$ with the semigroup $\zpn$. In the sequel, we note
$\{e_k, 1\leq k\leq N\}$ the natural basis of this semigroup.\par
For all $k$, define $E^{\tilde w_0}_{\beta_k}=E_{\beta_k}=T_{i_1}\ldots
T_{i_{k-1}}(E_{\alpha_{i_k}})$.
For all $\bfm=(m_i)\in\zpn$, set $E^{\tilde
w_0}(\bfm)=E(\bfm):=E_{\beta_1}^{(m_1)}\ldots
E_{\beta_N}^{(m_N)}$. It is known that $\{E(\bfm), \, \bfm\in\zpn\}$
is a basis of $\uqn$ called the Poincar\'e-Birkhoff-Witt basis, in short
PBW-basis,
associated to the reduced decomposition $\tilde w_0$. In the same way, we
can define the PBW-basis $\{F(\bfm), \, \bfm\in\zpn\}$ of $\uqnm$.
\par
In the sequel, we call a (left) factor of $\tilde w_0$ a reduced
decomposition
$\tilde w= s_{i_1}\ldots s_{i_{k}}$, $1\leq k\leq N$. Conversely, we say
that
$\tilde w_0$ is a (right) completion of $\tilde w$.\par
Let $w$ be in the Weyl group and let $\tilde w= s_{i_1}\ldots s_{i_{k}}$
be a reduced decomposition of $w$ which is
completed to $\tilde w_0=s_{i_1}s_{i_2}\ldots s_{i_N}$.
Let $U_q(\n_{\tilde w})$ be the $\Q(q)$-space generated by the
$E^{\tilde w_0}(\bfm)$ such that $m_i=0$ for $i>k$. By [7, 2.3],
$U_q(\n_{\tilde w})$ is a subalgebra of $\uqn$ which depends only on $w$ and
not
on the reduced decomposition $\tilde w$. In the sequel, we shall note it
simply $U_q(\n_{w})$.\par\noindent
Recall the following theorem, [6, 1.7]:
\begin{thm} Fix a reduced decomposition $\tilde w_0$ of $w_0$ and
set $${\cal F}_{\bfm}^{\tilde w_0}(\uqn)=\oplus_{\bfn\prec\bfm}\Q(q)
E^{\tilde w_0}({\bf n}),
\hskip 5mm
\bfm\in\zpn,$$
where $\prec$ is the right lexicographical ordering of $\zpn$.
Then, the spaces ${\cal F}_{\bfm}^{\tilde w_0}(\uqn)$,
$\bfm\in\zpn$, define a $\zpn$-filtration of $\uqn$.
The associated graded algebra Gr$^{\tilde w_0}(\uqn)$ is generated
by Gr$^{\tilde w_0}(E_\alpha)$, $\alpha\in\rp$ with relations (up to a
sign):
$$\hbox{Gr}^{\tilde w_0}(E_\alpha)\hbox{Gr}^{\tilde w_0}(E_\beta)
=q^{<\alpha,\beta>}\hbox{Gr}^{\tilde w_0}(E_\beta)\hbox{Gr}^{\tilde
w_0}(E_\alpha),
\hskip 5 mm \alpha<\beta.$$
\end{thm}
As in [11, 4.2], we define the bilinear forms $d^{\tilde w_0}=d$ and
$c^{\tilde w_0}=c$ on $\zpn\times\zpn$
such that
$$\hbox{Gr}(E(\bfm))\hbox{Gr}
(E(\bfn))
\in q^{-d(\bfm,\bfn)}(\pm 1+q\Z[q])\hbox{Gr}(E(\bfm+\bfn))$$
$$\hbox{Gr}(E(\bfm))\hbox{Gr}
(E(\bfn))=q^{c(\bfn,\bfm)}\hbox{Gr}
(E(\bfn))\hbox{Gr}(E(\bfm)),$$ up
to a sign. To be more precise, we have,
$$d(\bfm,\bfn)=\sum_{i>j}<\beta_i,\beta_j>m_in_j +
{1\over 2}\sum_i <\beta_i,\beta_i>m_in_i.$$
We define the $\Z[q]$-lattice $\lat$ generated by the
$\{E^{\tilde w_0}(\bfm), \, \bfm=(m_i)\in\zpn\}$. From [15, Proposition
2.3],
this lattice does not depend on the choice of $\tilde w_0$.
\sss  There exists,  see [23], a unique non degenerate Hopf pairing
$(\,,\,)$ on $\uqb\times\uqbm$ such that
$$(u^+,\,u_1^-u_2^-)=(\Delta(u^+),\,u_1^-\otimes u_2^-)\,,\hskip 1cm
u^+\in\uqb\,;\,u_1^-,u_2^-\in \uqbm$$
$$(u_1^+u_2^+,\,u^-)=(u_2^+\otimes u_1^+,\,\Delta(u^-))\,,\hskip 1cm
u^-\in\uqbm\,;\,u_1^+,u_2^+\in\uqb$$
$$(K_\lambda,\,K_\mu)=q^{-(\lambda,\mu)}\,,\hskip 1cm
\lambda,\mu\in Q$$
$$(K_\lambda,\,F_i)=0\,,\hskip 1cm\lambda\in Q\,,\,1\leq i\leq
n$$
$$(E_i,\,K_\lambda)=0\,,\hskip 1cm \lambda\in Q\,,\,1\leq i\leq
n$$
$$(E_i,\,F_i)={1\over 1-q_{\alpha_i}^{-2}}\delta_{ij}\,,\hskip 1cm 1\leq i,j
\leq
n,$$
where $\Delta$ is the comultiplication of the Hopf algebra $\uq$.\par
Let $\{E(\bfm)^*, \, \bfm\in\zpn\}$ be the dual basis of
$\{F(\bfm), \, \bfm\in\zpn\}$ for this pairing. \par
From [12], we claim that:
\begin{claim} For all $\bfm$ in $\zpn$ there exists a  function
$f_{\bfm}(q)$ in $\Q(q)$ such that $E(\bfm)^*=f_{\bfm}(q)E(\bfm)$
and $f_{\bfm}(0)=1$. Moreover,
this function is an eigenvector for the $\Q$-automorphism of $\Q(q)$ defined by $q\mapsto
q^{-1}$.
\end{claim}
We set $$\lat^*=\oplus_{\bfm\in\zpn}\Z[q]E(\bfm)^*.$$
\sss In this section, we give some results about Lusztig's canonical basis
$\B$ and its dual $\B^*$. \par\noindent
First of all, for $\lambda$ in $P^+$, let
$\lql$ be the simple $\uq$-module with highest weight $\lambda$. Choose
a highest weight vector $v_\lambda$. It is known that $\lql$
verifies the Weyl character formula. For all $w$ in $W$, let $v_{w\lambda}$
be a (non zero) extremal vector of weight $w\lambda$ and let
$V_{q,w}(\lambda):=\uqn.v_{w\lambda}$ be the Demazure module.\par
Let $\B$ be the Lusztig's canonical basis of $\uqnm$, [15], which
coincides with Kashiwara's global basis, [9]. It verifies the following
property:
\begin{thm} Fix $\lambda$ in $P^+$, and let $\B(\lambda):=\{b\in\B,\,
bv_\lambda\not=0\}$. Then, the set $\B(\lambda).v_\lambda$ is a basis of
$\lql$. Moreover, for all $w$ in $W$ there exists a unique subset $\B_w$
of $\B$, which does not depend on $\lambda$ and such that
$(\B(\lambda)\cap\B_w).v_\lambda$ is a basis of $V_{q,w}(\lambda)$.
\end{thm}
\begin{remark} In the sequel, if no confusion occurs, we shall identify
$\B(\lambda)$ with $\B(\lambda).v_\lambda$.
\end{remark}
Let $\B^*\subset\uqn$
be the dual basis in $\uqn$, i.e. $(b^*,b')=\delta_{b,b'}$. Remark that this
basis is not really canonical since it depends on the choice of a Hopf
pairing
$(\,,\,)$. Nethertheless, we shall call it the dual canonical basis.\par
Let $\eta$ be the $\Q$-automorphism of $\uq$ such that $\eta(E_i)=E_i$,
$\eta(F_i)=F_i$, and $\eta(q)=q^{-1}$.
Let $\sigma$ be the $\Q(q)$-antiautomorphism of
$\uq$ such that $\overline E_i=E_i$,
$\overline F_i=F_i$.
We can now give a characterization of $\B^*$, see  [11, Proposition 16].
This
characterization will give rise to the Lusztig's parametrization of the dual
canonical basis $\B^*$. It depends on the choice of a reduced decomposition
of $w_0$.
 \begin{proposition} Fix a reduced decomposition $\tilde w_0$ of $w_0$.
 Then, for each $\bfm$ in $\zpn$, there exists a unique homogeneous element
$X:=B^{\tilde w_0}(\bfm)^*=B(\bfm)^*$ in $\uqn$  such that
 $$\sigma\eta(X)=(-1)^{\tr(X)}q^{<\wt(X),\wt(X)>/2}q_{X}X,\;\;
 X\in E(\bfm)^*+q\lat^*,$$
 where $\wt(X)=\sum_ik_i\alpha_i$, is the weight of $X$,
$q_X=\prod_iq_{\alpha_i}^{k_i}$, and $\tr(X)=\sum_ik_i$.
\end{proposition}
\begin{remark} First remark that the proposition implies that $\lat$,
resp. $\lat^*$, is the $\Z[q]$-lattice generated by $\B$, resp. $\B^*$.
Remark also that the eigenvalue of $X=B(\bfm)^*$ for $\sigma\eta$ only
depends on the weight of $X$. Now, it can be easily seen that the first
condition
in the proposition can be replaced by ``$X$ is an eigenvector for
$\sigma\eta$".
\end{remark}
We note $B^{\tilde w_0}(\bfm)=B(\bfm)$ the corresponding element in the
canonical basis $\B$. For $w$ in $W$ and $\lambda$ in $P^+$, we set
$\B_w^*:=\{b^*\in \B^*, b\in\B_w\}$, and $\B(\lambda)^*:=
\{b^*\in \B^*, b\in\B(\lambda)\}$.
\section{Preliminary results on the
dual cano\-ni\-cal basis.}
\sss In this section, we obtain a property of triangularity for the decomposition
matrix of
the dual canonical basis in a general PBW basis.
\begin{proposition} Let $\tilde w_0$ be a reduced decomposition of $w_0$.
Then, for all $\alpha$ in $\rp$, the element $(E_\alpha^{\tilde w_0})^*$
of the dual PBW-basis belongs to $\B^*$.
\end{proposition}
\begin{proof}: By Proposition 1.4 and Remark 1.4, it is enough to
prove that $(E_\alpha^{\tilde w_0})^*$ is an eigenvector for $\sigma\eta$.
By Claim 1.3, it is sufficient to prove that this is true for
$E_\alpha^{\tilde w_0}$.\par
Set $\tilde w_0=s_{i_1}\ldots s_{i_N}$. We have
 $E_\alpha^{\tilde w_0}=T_{i_1}\ldots T_{i_{p-1}}(E_{i_p})$ for $p$ such
that $\alpha=\beta_p$.  By (1.2.1)-(1.2.3), we obtain that for all
homogeneous element $X$
in $\uqn$ of weight $\mu$, we have
$\sigma\eta T_i(X)=(-q_{\alpha_i})^{\tr(s_i\mu-\mu)}T_i(\sigma\eta(X)))$.
Remark that
$E_{i_p}$ is an homogeneous eigenvector for $\sigma\eta$. By induction, this
is
also true for $E_\alpha^{\tilde w_0}$.
\end{proof}
The following corollary generalizes [15, Theorem 9.13 (a)].
\begin{corollary} Let $\tilde w_0$ be a reduced decomposition of $w_0$.
Then, for all $\bfm$ in $\zpn$, we have the following homogeneous
decompositions.\par\noindent
\item{(i)} $\sigma\eta(E^{\tilde w_0}(\bfm))=e_{\bfm}^{\bfm}
E^{\tilde w_0}(\bfm)+\sum_{\bfn \prec\bfm}
e_{\bfm}^{\bfn}E^{\tilde w_0}(\bfn)$, $e_{\bfm}^{\bfn}\in
\Z[q,q^{-1}]$,\par\noindent
\item{(ii)} $B^{\tilde w_0}(\bfm)=F^{\tilde w_0}(\bfm)+\sum_{\bfm \prec\bfn}
d_{\bfm}^{\bfn}F^{\tilde w_0}(\bfn)$, $d_{\bfm}^{\bfn}\in
q\Z[q]$,\par\noindent
\item{(iii)} $B^{\tilde w_0}(\bfm)^*=E^{\tilde w_0}(\bfm)^*+\sum_{\bfn
\prec\bfm}
c_{\bfm}^{\bfn}E^{\tilde w_0}(\bfn)^*$, $c_{\bfm}^{\bfn}\in q\Z[q]$.
\end{corollary}
\begin{proof} Let $\bfm$ be in $\zpn$ and set $\sigma\eta(E(\bfm))
=\sum_{\bfn}e_{\bfm}^{\bfn} E(\bfn)$, $e_{\bfm}^{\bfn}\in\Z[q,q^{-1}]$.
Up to a multiplicative scalar, we have
from the previous proposition: $\sigma\eta(E(\bfm))=\
\sigma\eta(E_{\beta_1}^{(m_1)}
\ldots
E_{\beta_N}^{(m_N)})$
\par\noindent $=\sigma\eta(E_{\beta_N}^{(m_N)})\ldots
\sigma\eta(E_{\beta_1}^{(m_1)})=
E_{\beta_N}^{(m_N)}\ldots E_{\beta_1}^{(m_1)}$.
Hence, by Theorem 1.2,
Gr$^{\tilde w_0}(\sigma\eta(E(\bfm)))$\par\noindent
$=$Gr$^{\tilde w_0}(E(\bfm))$, up to a
multiplicative scalar. This implies (i).\par
By duality, it is sufficient to prove (iii). The coefficients
$c_{\bfm}^{\bfn}$ are in $q\Z[q]$ by Proposition 1.3. Now, the
property of triangularity comes from Proposition 1.4 and (i).
\end{proof}
\begin{remark} In the previous corollary, the lexicographical ordering
$\prec$ can be replaced by a coarser ordering: the one generated by
$\bfm\leq e_k+e_{k'},\; 1\leq k<k'\leq N\Leftrightarrow E(\bfm)$ is a term
of the PBW decomposition of
$E_{\beta_{k}}E_{\beta_{k'}}-q^{<\beta_{k'},\beta_{k}>}E_{\beta_{k'}}
E_{\beta_{k}}$.
\end{remark}
\sss The previous corollary implies the compatibility of the dual canonical
basis with the
space $U_q(\n_w)$, $w\in W$.
\begin{proposition} Let $w$ be in $W$. Then, $U_q(\n_w)$ is generated as a
space by a part of $\B^*$.
\end{proposition}
\begin{proof} Let $\tilde w$ be a reduced decomposition of $w$ and let
$\tilde w_0$ be a reduced decomposition of $w_0$ which completes $\tilde w$.
Fix an element of the PBW-basis which belongs to $U_q(\n_w)$. Then, by 1.2,
all smaller elements of the PBW-basis, with the same weight, belong to
$U_q(\n_w)$. By Corollary 2.1, this implies the proposition.
\end{proof}
\section{Quantum flag minors.}
\sss For all $\lambda$ in $P^+$, the weights of $\B(\lambda)^*$ are the
$\lambda-\mu$, where $\mu$ runs over the weights of $\lql$ (with
multiplicity). The quantum flag minors, see [10], [2, 4.2], are elements of
$\B(\varpi_i)^*$, $1\leq i\leq n$, which correspond to extremal vectors.
To be more precise, let $w$ be in $W$ and let $\tilde w=s_{i_1}\ldots s_{i_k}$ be a reduced
decomposition of $w$. There exists a unique element in $\B(\varpi_{i_k})^*$
with weight $(Id-w)(\varpi_{i_k})$. Note $\Delta_{\tilde w}^*$ this element. $\Delta_{\tilde w}^*$ is a quantum flag minor and each quantum flag
minor can be written in this way.\par
The following proposition generalizes some properties of the $q$-center of
$\uqn=U_q(\n_{w_0})$ proved in [4, Proposition 3.2] to $U_q(\n_{w})$.
\begin{proposition} Fix an element $w$ in $W$ and let $\tilde
w=s_{i_1}\ldots s_{i_k}$ be a reduced decomposition of $w$. We have
\par\noindent
\item{(i)} Let $X_\mu$ be an element of weight $\mu$ in $U_q(\n_w)$, then
$\Delta_{\tilde w}^*X_\mu=q^{<(Id+w)\varpi_{i_k},\mu>}X_\mu\Delta_{\tilde
w}^*$,
\par\noindent
\item{(ii)} let $b_\mu^*$ be an element of weight $\mu$ in
 $\B_w^*$, then
$q^{<\varpi_{i_k},\mu>}b_\mu^*\Delta_{\tilde w}^*\in\B^*$ mod $q\lat^*$.
\end{proposition}
\begin{proof}
Let's prove (i). Set $\tilde w'=s_{i_1}\ldots s_{i_{k'}}$, for $k'\leq k$.
From [5, (3.3.2)], we obtain
$$\Delta_{\tilde w}^*\Delta_{\tilde w'}^*=
q^{-<(Id-w)\varpi_{i_k}-2\varpi_{i_k},\mu'>}
\Delta_{\tilde w'}^*\Delta_{\tilde w}^*,$$
where $\mu':=(Id-w')\varpi_{i_{k'}}$ is the weight of $\Delta_{\tilde w'}^*$
and $k'<k$.
Hence
$$\Delta_{\tilde w}^*\Delta_{\tilde w'}^*=
q^{<(Id+w)\varpi_{i_k},\mu'>}
\Delta_{\tilde w'}^*\Delta_{\tilde w}^*.$$
This formula holds for $k'=k$ because
$<(Id+w)\varpi_{i_k},(Id-w)\varpi_{i_{k}}>=<\varpi_{i_k},\varpi_{i_k}>
-<w\varpi_{i_k},w\varpi_{i_k}>=0$, by $W$-invariance.
Now, (i) comes from the fact that the division ring generated by the
$\Delta_{\tilde w'}^*$, $1\leq k'\leq k$, is the division ring of
$U_q(\n_w)$, [5, corollary 3.2].\par
The proof of (ii) is a straightforward generalization of [4,
Proposition 3.2].
Let's sketch the proof. Let $b_\mu$ and $\Delta_{\tilde w}$ be the elements
in $\B$ corresponding respectively to $b_\mu^*$ and $\Delta_{\tilde w}^*$
and suppose $b_\mu\in \B(\lambda)$, $\lambda\in P^+$. Then,
$ \Delta_{\tilde w}
 \otimes b_\mu\in\B(\varpi_{i_k})\otimes\B(\lambda)$ and this element
corresponds
 to an element of the crystal basis at $q=0$.
 We know that $b_\mu$ is in $\B_w$. Hence, by Littelmann's
path model, we can associate to $b_\mu$ a chain of elements in $W$ which are
lower than $w$ for the Bruhat ordering. Moreover, the chain associated to
 $\Delta_{\tilde w}$ is reduced to $w$. So, both chains can be compared for
the Bruhat
 ordering. This implies that $ \Delta_{\tilde w}
 \otimes b_\mu\in\B(\varpi_{i_k}+\lambda)$ at the crystal level. This
implies
 (ii) by [11, Proposition 33].
 \end{proof}
 \sss By [7, Theorem 3.2] and Proposition 2.2, we have:
 \begin{lemma} Let $\tilde w_0= s_{i_1}\ldots s_{i_N}$
 be a reduced decomposition of $w_0$ and let
 $w=s_{i_1}\ldots s_{i_k}$, $1\leq k\leq N$. Let $\bfm$ be in
 $\zpn$ such that $m_{k'}=0$ for $k'>k$. Then, $B^{\tilde w_0}(\bfm)$ is in
$\B_w$.
 \end{lemma}
 Remark that this lemma can be directly proved by applying the
Berenstein-Zelevinsky
 formula, [2, Theorem 3.7], for the transition map between the Lusztig
parametrization
 and the string parametrization of the dual canonical basis. Indeed, by
[9, Theorem 12.4],
 the elements of $\B_w^*$ are characterized by a string parametrization.
 \par We can now prove the key proposition:
 \begin{proposition} Let $\tilde w_0= s_{i_1}\ldots s_{i_N}$
 be a reduced decomposition of $w_0$ and let
 $w=s_{i_1}\ldots s_{i_k}$, $1\leq k\leq N$. Fix $\bfm$  in
 $\zpn$. Then,
 $$q^{d(\bfnw,\bfm)}\Delta_{\tilde w}^*E^{\tilde w_0}(\bfm)^*\in
 E^{\tilde w_0}(\bfnw+\bfm)^*+q\lat^*,$$
  where $\bfnw$ is such that $\Delta_{\tilde w}^*=
 B^{\tilde w_0}(\bfnw)^*$.
 \end{proposition}
 \begin{proof}
 We first suppose that $\bfm$ is in $\Z_{\geq 0}^k\times\{0\}^{N-k}$. Then,
by
 Proposition 3.1 and the previous lemma, we have
 $$q^{<\nu-\varpi_{i_k},\mu>}\Delta_{\tilde w}^*B^{\tilde w_0}(\bfm)^*\in
 \B^*+q\lat^*,$$
 where $\nu$ and $\mu$ are respectively the weights of
 $\Delta_{\tilde w}^*$ and $B^{\tilde w_0}(\bfm)^*$.\par
Recall the notations of 1.2. By [11, 4.2], we have
$$d(\bfnw,\bfm)+d(\bfm,\bfnw)=<\nu,\mu>,\hskip 5mm
d(\bfnw,\bfm)-d(\bfm,\bfnw)=c(\bfnw,\bfm).$$
Using Proposition 3.1, we obtain
$d(\bfnw,\bfm)=<\nu-\varpi_{i_k},\mu>.$
Hence, there exists $\bfmp$ in $\zpn$ such that
$$q^{d(\bfnw,\bfm)}\Delta_{\tilde w}^*B^{\tilde w_0}(\bfm)^*
\in E^{\tilde w_0}(\bfmp)^*+\sum_{\bfn}q\Z[q]E^{\tilde w_0}(\bfn)^*.$$
Using Theorem 1.2 and Claim 1.3, we obtain $\bfmp=\bfnw+\bfm$.
By Proposition 2.2, this implies
$$q^{d(\bfnw,\bfm)}\Delta_{\tilde w}^*E^{\tilde w_0}(\bfm)^*\in
 E^{\tilde w_0}(\bfnw+\bfm)^*+q\lat^*.$$
 We can now study the general case. Let $\bfm$ be in $\zpn$ and decompose
 $\bfm=\bfn+\bfp$, with $\bfn\in\Z_{\geq 0}^k\times\{0\}^{N-k}$ and
 $p_l=0$ for $l\leq k$.
 We have $q^{d(\bfnw,\bfm)}\Delta_{\tilde w}^*E(\bfm)^*
 = q^{d(\bfnw,\bfn)}\Delta_{\tilde w}^*E(\bfn)^*E(\bfp)^*
\in (E(\bfnw+\bfn)^*+\sum q\Z[q]E(\bfr)^*)E(\bfp)^*,$
where $\bfr$ runs over $\Z_{\geq 0}^k\times\{0\}^{N-k}$.
Hence, $q^{d(\bfnw,\bfm)}\Delta_{\tilde w}^*E(\bfm)^*
\in E(\bfnw+\bfn+\bfp)^*+\sum q\Z[q]E(\bfr+\bfp)^*
\subset E(\bfnw+\bfm)^*+q\lat^*.$
This ends the proof.
\end{proof}
\section{Quiver orientations and quantum flag minors}
According to [4], Proposition 3.2 is almost what we need if we
want to test the Berenstein-Zelevinsky
conjecture when one element is a quantum flag minor. In fact, we have
to prove an analogue of this proposition where $E^{\tilde w_0}(\bfm)^*$
is replaced by $B^{\tilde w_0}(\bfm)^*$. This can be obtained if we prove
some increasing
property of the linear form $d^{\tilde w_0}(\bfnw,?)$. By results of
Reineke, [19], and Bongartz, [3], it is possible to
prove those properties by the quiver approach, [20], of quantum groups.
\sss This section refers to [21] for notations and definitions.  Define
the graph $\ga$ as follows: if $\g$ is simply laced, i.e.
of type A-D-E, $\ga$ is the Coxeter graph of $\g$, if $\g$ is not simply
laced then $\ga$ is the A-D-E graph such that the graph of $\g$ is a
quotient of this graph, see [16, par. 14], [18, par. 7]. Fix an orientation
$\gaf$
of $\ga$ and let Mod$k\gaf$ be the category of $k$-representations of
$\gaf$, where $k$ is an algebraically closed field.\par
For $M$, $N$ in Mod$k\gaf$, set 
$$\varepsilon(M,N)=\dim_k\Hom(M,N),\hskip 5mm
 \zeta(M,N)=\dim_k\Ext^1(M,N).$$ Let Ind$k\gaf$ be the set of
indecomposable modules of Mod$k\gaf$. Let $\tau$: Ind$k\gaf\rightarrow$
Ind$k\gaf\cup\{0\}$ be the Auslander-Reiten translation, with the convention
that $\tau(M)=0$ if $M$ is a projective module. Recall, [21], that
$$\zeta(M,N)=\varepsilon(N,\tau(M)), \hskip 5mm
 M,N \in
\hbox{Ind}k\gaf.\leqno (4.1.1)$$
By the work of Ringel, [20], the algebra $\uqn$ can be realized as a
deformation of the Hall algebra associated to $\gaf$. In particular, for a
special reduced decomposition $\tilde w_0(\gaf)$ of $w_0$, the PBW-basis of
$\uqn$ can be naturally parametrized by Mod$k\gaf$. The following
proposition is a recollection of results which can be found in [19],
[15, par 4]. It gives informations on the link  between the
Mod$k\gaf$-parametrization and the Lusztig parametrization associated to
$\tilde w_0(\gaf)$.
\begin{proposition} There exists a reduced decomposition $\tilde
w_0(\gaf)=s_{i_1}s_{i_2}\ldots s_{i_N}$, defined up to commutations, and a
$\Z_{\geq 0}$-linear bijection $\iota$:  Mod$k\gaf\rightarrow\zpn$ such
that:\par\noindent
\item{(i)} if $M$ is indecomposable non projective and if
$\iota(M)=e_k$ , then $\iota(\tau(M))=e_{k'}$,
where $k'<k$ is maximal such that $i_{k'}=i_k$,
\item{(ii)}
$d(\iota(M),\iota(N))=\varepsilon(N,M)-\zeta(M,N)$.\par\noindent
Moreover, let $i$ be a sink of $\gaf$ then, up to commutations, $i_1=i$. Let
$\gaf_i$ be the oriented graph obtained from $\gaf$ by reversing the arrows
which go to $i$, we have $\tilde w_0(\gaf_i)=s_{i_2}\ldots s_{i_N}s_{i^*}$,
where $i^*$ is such that $w_0(\alpha_i)=-\alpha_{i^*}$.
\end{proposition}
\sss In the following sections, we shall replace the lexicographical
ordering $\prec$ by the coarser ordering discussed in Remark 2.1. We still
note it $\prec$. If the reduced decomposition $\tilde w_0$ is associated to
a graph orientation, then this ordering is the Ext-ordering on
Mod$k\gaf\simeq\zpn$, [3], which is also the degeneration ordering,
[{\it loc. cit.}, Corollary 4.2,  Proposition 3.2]. We are now ready to
prove the expected ``increasing
property" of $d^{\tilde w_0}(\bfnw,?)$.
\begin{proposition} Let $\tilde w_0(\gaf)=s_{i_1}s_{i_2}\ldots s_{i_N}$
be the reduced decomposition of $w_0$ associated to $\gaf$ and let
$\tilde w=s_{i_1}s_{i_2}\ldots s_{i_k}$ be a subword of $\tilde w_0(\gaf)$.
Let $\bfnw$ be as in 3.1. Then, the form $d(\bfnw,?)$ is increasing for
$\prec$.
\end{proposition}
\begin{proof} By [5, 3.2.2], we have $\bfnw=\sum e_l$ where $l$ runs over
the
set of integers of $[1,k]$ such that $i_l=i_k$. From Proposition 4.1, we
obtain that $$d(\bfnw,?)=\sum_{i\geq 0}\varepsilon(?,\tau^i(M_k))
-\sum_{i\geq 0}\zeta(\tau^i(M_k),?),$$ where $M_k$ is the indecomposable
module corresponding to $e_k$. By (4.1.1), we obtain after elimination that
$d(\bfnw,?)=\varepsilon(?,M_k)$, which is increasing by [3,
Proposition 3.2].
\end{proof}
\sss Now, the natural question is to know if all quantum flag minor can be
associated to a graph decomposition $\gaf$. In clear, fix  $u$ in $W$ and a
reduced decomposition $\tilde u$ of $u$. Is there an orientation $\gaf$ of
$\ga$ and a subword $\tilde w$ of $\tilde w_0=\tilde w_0(\gaf)$ such that
$\Delta_{\tilde u}^*=\Delta_{\tilde w}^*$. This is true for the A$_n$
case, but not in general:
\begin{claim} Let $\g$ be of type A$_n$. Then, all quantum flag minors can
be realized as $B^{\tilde w_0}(\bfnw)^*$ where $\tilde w_0$ is associated to
an orientation of $\ga$.
\end{claim}
\begin{proof}
It is known that for all $k$, $1\leq k\leq n$, the set of quantum minors in
$\B(\varpi_k)^*$ is naturally indexed by a set of lines
$I=\{i_1<i_2<\ldots<i_k\}$. We generally note it $\Delta_q(I,J)$ where
$J=\{1,\ldots,k\}$ is the set of its columns. This (flag) quantum minor
corresponds to an extremal vector in $V_q(\varpi_k)$ associated to
$$w=s_{(i_1-1)}\ldots s_1s_{(i_2-1)}\ldots s_1\ldots
s_{(i_k-1)}\ldots s_1
W^{k},$$
where $W^{k}:=\{u\in W,\,u\varpi_{k}=\varpi_{k}\}$.
\par\noindent
We claim that this reduced decomposition 
$\tilde w=s_{(i_1-1)}\ldots s_1s_{(i_2-1)}\ldots s_1\ldots
s_{(i_k-1)}\ldots s_1$ can be completed to a reduced
decomposition of $w_0$ which is associated to an orientation of $\ga$.
Indeed, this follows from the last assertion of Proposition 4.1 and the well
known fact: $\tilde w_0:=s_1s_2s_1s_3s_2s_1\ldots s_ns_{n-1}\ldots s_2 s_1$
is a reduced decomposition associated to the orientation $\gaf$ of A$_n$ such that
all arrows are oriented to the left. Indeed, with the notations of Proposition 4.1, $\tilde w$ can be completed to the reduced decomposition of $w_0$ associated to the quiver $\gaf_{{\bf i}}$, where ${\bf i}$ is the following sequence of reversing arrows
$${\bf i}=121\ldots(n-k)\ldots 1(n-k+1)\ldots i_1(n-k+2)\ldots i_2\ldots n\ldots i_k.$$
\end{proof}
\begin{remark} The claim is not true if we take $\g$ of type D$_4$. With the
standard
notations, the quantum flag minor associated to the reduced decomposition
$\tilde w=s_2s_1s_3s_2$ can not be realized from a quiver orientation.
\end{remark}
\section{Adapted algebras associated to a quiver orientation.}
\sss We first recall some basic facts on adapted algebras associated to a
reduced decomposition of $w_0$, see [4]. Fix a reduced decomposition
$\tilde w_0$ of $w_0$ and let $\tilde W$ be the set of left subwords of
$\tilde w_0$. Then, the quantum flag minors $\Delta_{\tilde w}^*$,
$\tilde w\in \tilde W$, $q$-commute. Moreover, they generate a
$\Q(q)$-algebra
$A_{\tilde w_0}$ such that\par\noindent
1) $A_{\tilde w_0}$ is a $q$-polynomial algebra of GK-dimension
$N=\hbox{ Card }\tilde W$,\par\noindent
2) As a space, $A_{\tilde w_0}$ is generated by a part of the dual canonical
basis, namely the monomials in the $\Delta_{\tilde w}^*$,
$\tilde w\in \tilde W$.\par
Now, we reach the main theorem:
\begin{thm} Fix an orientation $\gaf$ of $\ga$ and set $\tilde w_0=\tilde
w_0(\gaf)$. Suppose that $b^*$  in $A_{\tilde w_0}\cup\B^*$ and $(b')^*$ in
$\B^*$ $q$-commute, then they are multiplicative, i.e. their product is an
element of the dual canonical basis up to a power of $q$.
\end{thm}
\begin{proof} Set $\tilde w_0=s_{i_1}s_{i_2}\ldots s_{i_N}$. The elements
$b^*=B^{\tilde w_0}(\bfm)^*$ and $(b')^*=B^{\tilde w_0}(\bfmp)^*$ $q$-commute.
Hence, in order to prove that they are multiplicative, it is sufficient, see
[11, 5.1], to prove that
$$ q^{d^{\tilde w_0}(\bfm,\bfmp)}B^{\tilde w_0}(\bfm)^*B^{\tilde
w_0}(\bfmp)^*\in B^{\tilde w_0}(\bfm+\bfmp)^*+q{\cal L}^*.\leqno (5.1.1)$$
Recall that, by Corollary 2.2, $$B^{\tilde w_0}(\bfmp)^*=E^{\tilde
w_0}(\bfmp)+\sum_{\bfn\prec\bfmp}c_{\bfmp}^{\bfn}E^{\tilde
w_0}(\bfn)^*,\leqno (5.1.2)$$
with $c_{\bfmp}^{\bfn}\in q\Z[q]$. Moreover, $B^{\tilde w_0}(\bfm)^*$ is a
product of $\Delta_{\tilde w}^*$,
where  $\tilde w$ runs over a subset $\tilde W_0$ of $\tilde W$. This
implies that $d^{\tilde w_0}(\bfm,?)=\sum d^{\tilde w_0}(\bfnw,?)$, where
$\tilde w$ runs over $\tilde W_0$ with multiplicity. By Proposition 4.2,
this asserts that the form $d^{\tilde w_0}(\bfm,?)$ is increasing for
$\prec$. Hence, (5.1.1) is a consequence of Proposition 3.2 and (5.1.2).
\end{proof}
\begin{remark} It is likely that the theorem works for all reduced
decomposition
of $w_0$. In fact, all we have to prove is a generalization of Proposition 4.2 for general $\tilde w_0$.
\end{remark}
By Claim 4.3, we deduce the following corollary which generalizes a result
of [11]:
\begin{corollary} Let $\g$ be of type A$_n$. Let $b^*$ and $(b')^*$ be
$q$-commuting elements in $\B^*$. If $b^*$ is a quantum flag minor, then
$b^*$ and $(b')^*$ are multiplicative.
\end{corollary}
\sss
In order to understand a ``multiplicative" description of the dual canonical basis, S. Fomin and A. Zelevinsky have defined the so-called cluster algebras, [8]. Roughly speaking, cluster algebras are algebras equipped with a distinguished set of generators called cluster variables and this set is divided into a union of subsets called clusters. The cluster variables verify the so-called exchange relations. For each symmetrizable Cartan datum, Fomin and Zelevinsky associate a cluster algebra.\par
As a conclusion, we would like to enlight the link between the results of 5.1 and the theory of cluster algebras. This follows some discussions with A. Zelevinsky.  One of the most amazing fact in the Fomin-Zelevinsky theory is that many algebras, as $\C[U^-]$, $\C[G]$, encountered in the representation theory of semi-simple Lie groups seem to be realized in the algebraic framework of cluster algebras. In particular the subalgebras $A_{\tilde w_0}$ should specialize at $q=1$ onto a subalgebra of $\C[U^-]$ generated by a cluster. The connection with our problem is the following : it is reasonable to think  that if $b$ and $b'$ are $q$-commuting elements of the dual canonical basis, and if $b$ is specialized at $q=1$ onto a cluster variable, then $b$ and $b'$ are multiplicative. What we proved in 5.1 is a particular case of this conjecture when the cluster is associated to some reduced word.

We would like to thank Bernard Leclerc, Eric Vasserot and Andrei Zelevinsky for usefull conversations.

\end{document}